\documentclass[reqno,showkeys]{amsart}
\usepackage{amsfonts}
\usepackage{amssymb}
\usepackage{amsthm}
\usepackage{upref}
\makeatletter
\def\LaTeX{\leavevmode L\raise.42ex
    \hbox{\kern-.3em\size{\sf@size}{0pt}\selectfont A}\kern-.15em\TeX}
 \makeatother

\newcommand{\HH}{\textsf{H}}
\newcommand{\supp}{\operatorname{supp}}
\newcommand{\e}{\eqref}

\renewcommand\Im{\operatorname{Im}}

\numberwithin{equation}{section}

\newtheorem{lemma}{Lemma}[section]
\newtheorem{theorem}[lemma]{Theorem} 
\newtheorem{corollary}[lemma]{Corollary}
\newtheorem{proposition}[lemma]{Proposition}

\theoremstyle{definition}

\theoremstyle{remark}

\renewcommand{\det}{\operatorname{Det}}
\newcommand{\tr}{\operatorname{Tr}}

\def\qqq{\mathrel{\subset\mkern-15mu\lower.38ex\hbox{${\scriptscriptstyle\rightarrow}$}}}

\let\goth\mathfrak
\let\cal\mathcal

\let\Bbb\mathbb

\begin{document}
\title{A trace formula  for the Dirac operator}
\author{ D. R. Yafaev}
\address{ IRMAR, Universit\'{e} de Rennes I\\ Campus de
  Beaulieu, 35042 Rennes Cedex}
\email{yafaev@univ-rennes1.fr}
\begin{abstract}
Our goal is to extend the theory of the spectral shift function to the case where only the difference of  some powers of the resolvents  of self-adjoint operators belongs to the trace class. As an example, 
 we consider a couple of Dirac operators.
\end{abstract}
\maketitle

\thispagestyle{empty}

\section{Introduction}
The concept of the spectral shift function (SSF) first appeared in the work of I. M. Lifshits
\cite{Lif} in connection with the quantum theory of crystals. A mathematical theory of the SSF was shortly constructed by 
  M. G. Kre\u{\i}n in \cite{Kr1}. One of his results can be formulated in the following way.
Let $ H_0$ and $H$ be self-adjoint  operators with a trace class difference $V=H-H_0$ . Then there exists the function
$\xi(\lambda)=\xi(\lambda;H,H_0)$, $\xi\in L_1({\Bbb R})$, known as
the spectral shift function (SSF)    such that the trace formula
\begin{equation}
  \tr\Bigl( f(H) - f (H_0)\Bigr)=\int_{-\infty}^\infty \xi(\lambda) f^\prime(\lambda)
  d\lambda,\quad \xi(\lambda)=\xi(\lambda;H,H_0),
\label{eq:TF}\end{equation}
holds at least for all functions $f \in C_0^\infty({\Bbb R})$. Later in \cite{Kr2}
 M. G. Kre\u{\i}n  returned to this problem. In particular, he has shown that formula
similar to  (\ref{eq:TF}) remains true for a couple of unitary operators $U_{0}$ and $U$ with  a	trace class difference. In terms of self-adjoint operators this means that formula
 (\ref{eq:TF}) holds if the difference of the resolvents of the operators $H_{0}$ and $H$ belongs to the trace class $\mathfrak{S}_{1}$ (such operators
 $H_{0}$ and $H$ are called  resolvent comparable).
 This allows to write formula (\ref{eq:TF}) for sufficiently large class of differential operators. A relatively detailed presentation of the theory of the  SSF can be found in 
 \cite{BY2} or \cite{I}.

 In 1957 \cite{Ka3,Ka4,Ros} T. Kato and M. Rosenblum have proven the existence
  of  the wave operators (all relevant definitions can be found in \cite{I}) for the pair of self-adjoint operators with a trace class difference. This result was extended in  \cite{BK} by M. Sh. Birman and M. G. Kre\u{\i}n
  to unitary operators. More important, in the same paper they have found a connection
  between the scattering theory and the theory of the SSF. Actually, they have shown that   the corresponding scattering matrix $ S(\lambda;H,H_{0})$
  minus the identity operator $I$ belongs to  the trace class
   and
  \begin{equation}
\det S(\lambda;H,H_{0})=e^{-2\pi i\xi(\lambda;H,H_{0})}
\label{eq:Sxi}\end{equation}
for almost all $\lambda$ from the core of the spectrum of the operator $H_0$.

Later,  in \cite{Ka5} T. Kato has proven the existence of the wave operators for the pair $H_{0}$, $H$
under the assumption that
\begin{equation}
R^m(z)-R_0^m(z)\in \mathfrak{S}_{1},\quad \Im z \neq 0.
\label{eq:Rm}\end{equation}
 However the construction of the  SSF under this assumption seemed to be an open problem. Our goal here is to fill in this gap. Apart from its conceptual naturalness, there are very simple applications which require such an extension of the theory. 
 Here we give only one application -- to the Dirac operator.

We note that very general conditions of the existence of wave operators were obtained by  M. Sh. Birman  in the framework of the local trace class approach \cite{B2}. On the other hand, the local theory of the SSF due to L.S. Koplienko \cite{Kopl} is somewhat less satisfactory.

The author  is grateful to M. Z. Solomyak for useful discussions of the Double Operator Integrals theory.

\section{The spectral shift function}
Let us first formulate the result of M. G. Kre\u{\i}n
 obtained in \cite {Kr2}.
 
\begin{theorem}\label{KR}
Let 
\begin{equation}
(h-i)^{-1} - (h_{0}-i)^{-1}\in \mathfrak{S}_{1}.
\label{eq:R1}\end{equation}
Suppose that a function $g(\mu)$   has two bounded derivatives and
 $$
  \partial^\alpha( g(\mu)- g_{0}\mu^{-1})=O(|\mu |^{-1-\epsilon-\alpha}),
\quad \alpha=0,1,2, \quad \epsilon>0,
$$
where the constant $g_{0}$ is the same for $\mu\rightarrow \infty$ and
$\mu\rightarrow -\infty$. 
Then 
\[
g(H)- g(H_{0})\in \mathfrak{S}_{1} 
\]
and there exists the SSF $\xi(\mu;h,h_{0})$ such that
\begin{equation}
\int_{-\infty}^\infty |\xi(\mu;h,h_{0})|(1+|\mu|)^{-2}d\mu<\infty 
\label{eq:xi1}\end{equation}
and
\begin{equation}
  \tr\Bigl( g(h)- g(h_{0})\Bigr)
=\int_{-\infty}^\infty \xi(\mu;h,h_0)  g^\prime(\mu)
 d\mu.
\label{eq:TFG}\end{equation}
\end{theorem}
 
Our goal is to extend this result to the case where condition
\eqref{eq:Rm}  is satisfied for some (not necessarily $m=1$) odd $m$.

  \begin{theorem}\label{KRY}
 Let, for a pair of self-adjoint operators $H_{0}$ and $H$,
  the assumption \eqref{eq:Rm}  hold for some odd $m$ and all $\Im z \neq 0$.
Let a  function $f(\lambda)$   have two bounded derivatives, and let
$$
   \partial^\alpha( f(\lambda)- f_{0}\lambda^{-m})=O(|\lambda |^{-m-\epsilon-\alpha}),
\quad \alpha=0,1,2, \quad \epsilon>0,
$$
where the constant $f_{0}$ is the same for $\lambda\rightarrow \infty$ and
$\lambda\rightarrow -\infty$.
 Then the inclusion
 \begin{equation}
f(H)-f(H_{0})\in \mathfrak{S}_{1} 
\label{eq:Rmf}\end{equation}
  holds and there exists a function $($the SSF$)$
 $\xi(\lambda;H,H_{0})$ satisfying condition 
 \begin{equation}
\int_{-\infty}^\infty |\xi(\lambda;H,H_{0})|(1+|\lambda|)^{-m-1}d\lambda<\infty  
\label{eq:xi2}\end{equation} 
  such that the trace formula \e{eq:TF}  is true.
  Moreover, for the  corresponding scattering matrix $S(\lambda;H,H_{0})$, the operator $S(\lambda;H,H_{0})-I\in\mathfrak{S}_{1}$ and  relation \e{eq:Sxi}
 holds
for almost all $\lambda$ from the core of the spectrum of the operator $H_0$.
 \end{theorem}

Our proof of Theorem~\ref{KRY} relies on its
 reduction to Theorem~\ref{KR} although, similarly to \cite {Kr2}, we could have deduced it from the corresponding result for unitary operators.  Actually,  we construct   a function $ \varphi$ such that the operators $h_{0}=\varphi(H_{0})$ and $h =\varphi(H )$ are resolvent comparable. 
 To be more precise, we shall prove the following result.

\begin{theorem}\label{Rm}
Let the assumption \eqref{eq:Rm}  hold for some odd $m$ and all $\Im z \neq 0$.
Let $\varphi\in C^2({\Bbb R})$ and let $\varphi(\lambda)=\lambda^m$ for sufficiently large $|\lambda|$. Then the pair $h_{0}=\varphi(H_{0})$, $h =\varphi(H )$
satisfies condition \e{eq:R1}.
\end{theorem}

We postpone the proof of Theorem~\ref{Rm} until the next  section. Here we use it for the construction of the SSF.  Since, under the assumptions of  Theorem~\ref{KRY}, the operators  $h_{0}=\varphi(H_{0})$ and $h =\varphi(H )$ are resolvent comparable, we can apply Theorem~\ref{KR}
to the couple $h_{0}$, $h$.
Then we   define the SSF for the pair $H_{0}$, $H$ by the relation 
\begin{equation}
 \xi(\lambda;H,H_{0})= \xi( \varphi(\lambda);\varphi(H),\varphi(H_{0})).
\label{eq:SSFi}\end{equation}
Suppose that the function $\varphi$ is invertible  and that
\[
\varphi^\prime(\lambda)\geq c>0.
\]
Set $\mu=\varphi(\lambda)$,  $\psi=\varphi^{-1}$,  $g(\mu)=f(\psi(\mu))$.
Then formula (\ref{eq:TFG}) implies that
$$
  \tr\Bigl( f(H)- f(H_{0})\Bigr)
=  \tr\Bigl( g(h)-g(h_{0})\Bigr)
=\int_{-\infty}^\infty  
\xi(\mu;h,h_{0}) g^\prime (\mu) d\mu.
$$
This coincides with (\ref{eq:TF}) if the SSF $\xi(\lambda;H,H_{0})$ is defined by formula
(\ref{eq:SSFi}).

It follows from estimate   (\ref{eq:xi1}) and the conditions on $\varphi(\lambda)$ that the function \eqref{eq:SSFi} satisfies estimate
\eqref{eq:xi2}. The class 
  of functions $f(\lambda)$ for which formula (\ref{eq:TF})
is true is obtained from the class 
  of functions $g(\mu)$ by the change of variables $\mu=\varphi(\lambda)$.
  Finally, formula \eqref{eq:Sxi} follows from the same formula for the pair $h_{0}$, $h$, definition \eqref{eq:SSFi} of the SSF and the invariance principle for scattering matrices.
  This concludes the proof of Theorem~\ref{KRY} given Theorem~\ref{Rm}.  

\medskip

We emphasize that the trace formula \eqref{eq:TF} fixes the SSF
$ \xi(\lambda;H,H_{0})$ up to an additive constant only. This constant remains undetermined also by condition \eqref{eq:xi2}. Nevertheless (see \cite{I}) defining
the SSF via the corresponding perturbation determinant fixes $ \xi(\lambda;H,H_{0})$ up to an {\em integer} constant. It is exactly for this  ``correct" choice of a constant that relation \e{eq:Sxi}
 holds.

\section{Proof of Theorem~\ref{Rm}}
 Our proof of Theorem~\ref{Rm} relies on the theory of Double Integral Operators (DOI).
 Let  us briefly recall its basic notions (see \cite{BS1,BS2}, for details).
 Let $H_{0}$ and $H$ be a pair of self-adjoint operators in a Hilbert space ${\cal H}$.
 Denote by $E_{0} $ and $E$ their spectral families. Let us use the following result which can be deduced from \cite{BS1,BS2}.

\begin{proposition}\label{BS}
 Let
 \begin{equation}
K=\Phi (T)=\int_{-\infty}^\infty \int_{-\infty}^\infty
 K(\lambda,\mu)
dE(\mu) T dE_{0}(\lambda),
\label{eq:DOI1}\end{equation}
where the kernel $ K(\lambda,\mu)$ is bounded, i.e.,
   \begin{equation}
| K(\lambda,\mu)|  \leq C <\infty,
\label{eq:DOI2}\end{equation}
it is differentiable in $\lambda$ and
  \begin{equation}
  |\partial K (\lambda,\mu)/\partial \lambda| \leq C  (1+\lambda^2)^{-1}.
\label{eq:DOI2x}\end{equation}
Assume, moreover, that
 \begin{equation}
\lim_{\lambda\rightarrow+\infty}  K(\lambda,\mu) 
 =\lim_{\lambda\rightarrow-\infty}  K(\lambda,\mu) 
\label{eq:DOI3}\end{equation}
$($these limits exist by virtue of \eqref{eq:DOI2x}$)$.
Then the transformer $\Phi:\mathfrak{S}_{1}\rightarrow\mathfrak{S}_{1}$
defined by \e{eq:DOI1} is bounded.
  \end{proposition}
  
  Let us set
 \[
 g(\lambda)= g_{z}(\lambda)=(\lambda-z)^{-m}
 \]
 and
 \[
 T= T_{z}= g_{z}(H)- g_{z}(H_{0}).
 \]
 We need the   representation of the difference $f(H)-f(H_{0})$ in terms of the DOI:
\begin{equation}
f(H)-f(H_{0})=\int_{-\infty}^\infty \int_{-\infty}^\infty\frac{f(\lambda)-f(\mu)}{g(\lambda)-g(\mu)}
dE(\mu) T dE_{0}(\lambda).
\label{eq:DOI}\end{equation}
A technical problem with a proof of Theorem~\ref{Rm} is that the denominator
$g(\lambda)-g(\mu)$ in \e{eq:DOI} has  ``extra" zeros on the antidiagonal $\lambda=-\mu$. Therefore we shall split Theorem~\ref{Rm} into two separate assertions.

\begin{proposition}\label{Rm1}
Let $f\in C_{0}^2({\Bbb R})$, and let $\supp f\subset [-r,r]$. 
Let the assumption \eqref{eq:Rm}  be satisfied for some odd $m$ and  $  z =ia$  with a sufficiently large $($compared to $r)$ $a$. Then inclusion \eqref{eq:Rmf} holds.
\end{proposition}

\begin{proposition}\label{Rm2}
Let  $ \theta\in C ^2({\Bbb R})$,  $ \theta (\lambda)=0$ for $|\lambda|\leq r$ and
$ \theta (\lambda)=1$ for $|\lambda|\geq 2r$. Set
  \begin{equation}
f(\lambda)=\theta (\lambda) (\lambda^m-i)^{-1}.
\label{eq:DOI4}\end{equation}
 Then inclusion \eqref{eq:Rmf} holds if
 the assumption \eqref{eq:Rm}  is satisfied for some odd $m$ and  $  z =ia$  with a sufficiently small $($compared to $r)$   $a$. 
\end{proposition}

Propositions~\ref{Rm1} and \ref{Rm2} imply of course  Theorem~\ref{Rm}. Indeed,
 \begin{equation}
 (\varphi(H)-i)^{-1} -  (\varphi(H_{0})-i)^{-1}
 =(f_{0}(H)- f_{0}(H_{0}))+ (f(H)- f(H_{0})),
\label{eq:DOI4bi}\end{equation}
where $f_{0}(\lambda)=(1-\theta (\lambda)) (\varphi(\lambda)-i)^{-1}$ has finite support and $f$ is given by formula \e{eq:DOI4}. Both terms in the right-hand side of 
\e{eq:DOI4bi} belong to the trace class.

For the proof of Propositions~\ref{Rm1} and \ref{Rm2} we need the following elementary result.

\begin{lemma}\label{Rmx}
Set
 \begin{equation}
p(\lambda,\mu;z)=  (\lambda- z)^{m-1}+  (\lambda-z)^{m-2}(\mu-z)+\cdots+
 (\mu-z)^{m-1}.
\label{eq:DOI5}\end{equation}
Let $m$ be odd, and let $r$ be some fixed number.
Then for $|\lambda|\leq r$,  $|\mu|\leq r$, $z=ia$ and a sufficiently large $a$
\begin{equation}
|p(\lambda,\mu;z)|\geq c>0.
\label{eq:DOI6}\end{equation}
Similarly,  if either $|\lambda|\geq r$ or  $|\mu|\geq r$, $z=ia$  and $a$ is sufficiently small,  then  
\begin{equation}
|p(\lambda,\mu;z)|\geq c (|\lambda|+|\mu|)^{m-1}, \quad c>0.
\label{eq:DOI6a}\end{equation}
\end{lemma}

 \begin{proof}
By virtue of the equality 
$$
p(\lambda,\mu;z)=  (\lambda-ia)^{m-1}(1+\sigma+\cdots +\sigma^{m-1}),
 \quad \sigma =(\mu-ia)(\lambda-ia)^{-1},
 $$
both estimates (\ref{eq:DOI6}) and (\ref{eq:DOI6a})
(by the proof of  (\ref{eq:DOI6a}) we assume that $|\mu|\leqÊ|\lambda|$) reduce to the same estimate
 \begin{equation}
|1+\sigma+\cdots +\sigma^{m-1}| \geq c>0.
\label{eq:DOI8}\end{equation}
If  $|\lambda|\leq r$,  $|\mu|\leq r$ and  $a$ is  sufficiently large, then
\[
 \sigma =(1+i\mu a^{-1})(1+i \lambda a^{-1})^{-1}
 \]
 belongs to a neighbourhood of the point $1$ which implies (\ref{eq:DOI8}).
 
 Remark that the zeros of the function
 $1+\sigma+\cdots +\sigma^{m-1}$
 are given by the formula 
 $$\sigma_{k}=\exp(2\pi i k/m), \quad k=1,\ldots, m-1.$$
 If  $|\lambda|\geq r$,  $|\lambda|\geq |\mu| $ and  $a$ is  sufficiently small, then 
 the values of 
$$ \sigma=(x-i\varepsilon)(1-i \varepsilon)^{-1}, \quad
x=\mu\lambda^{-1}\in[-1,1], \quad \varepsilon=a \lambda^{-1},
$$
   belong  to a neighbourhood of the real axis and therefore
    are separated from all points $\sigma_{k}$.
   This again implies (\ref{eq:DOI8}).
  \end{proof}

  For the proofs of Propositions~\ref{Rm1} and \ref{Rm2},
  we consider DOI (\ref{eq:DOI}) with kernel
    \begin{equation}
K(\lambda,\mu)=   \frac{f(\lambda)-f(\mu)}{g(\lambda)-g(\mu)} 
\label{eq:DO}\end{equation}
and verify the assumptions of Propositions~\ref{BS}.
Obviously, the function $K(\lambda,\mu)$ tends to $f(\mu)g(\mu)^{-1}$ as $\lambda\rightarrow\pm\infty$
so that condition \e{eq:DOI3} is satisfied.
The estimates \eqref{eq:DOI2},  \eqref{eq:DOI2x}  we shall check separately 
under the assumptions of Propositions~\ref{Rm1} and \ref{Rm2}. Note previously that,
under the assumptions of Proposition~\ref{Rm1},
$K (\lambda,\mu)=0$  if  $|\lambda|\geq r$ and 
 $|\mu|\geq r$. Similarly, under the assumptions of Proposition~\ref{Rm2},
 $K (\lambda,\mu)=0$ 
if $|\lambda| \leq r$ and $|\mu|\leq r$.

Since
   \begin{equation}
g(\lambda) -g(\mu)=    (\lambda-z)^{-m}  (\mu-z)^{-m}(\mu-\lambda ) p(\lambda,\mu;z),
\label{eq:DO1}\end{equation}
we have that
   \begin{equation}
K(\lambda,\mu)=  - \Phi (\lambda,\mu) G(\lambda,\mu),
\label{eq:DO2}\end{equation}
where
  \begin{equation}
\Phi (\lambda,\mu)=   \frac{f(\lambda)-f(\mu)} { \lambda-\mu } 
  \label{eq:DO3}\end{equation}
  and
   \begin{equation}
G(\lambda,\mu)=    
    \frac{  (\lambda-z)^{m}  (\mu-z)^{m}}{ p(\lambda,\mu;z) }.
\label{eq:DO4}\end{equation}

  Let us start with Proposition~\ref{Rm1} and consider first the region
where $|\lambda|\leq R$, 
 $|\mu|\leq R$ for some $R$. Here we proceed from representation \e{eq:DO2}.
  The function  \e{eq:DO3} 
 is bounded because $f \in C^1$.
 Since 
  $$
| f(\mu)-f(\lambda) -f^\prime(\lambda)(\mu-\lambda ) |\leq 
2^{-1}(\mu-\lambda)^2 \supp_{|\nu|\leq R }|f^{\prime \prime}(\nu)|,
$$
 the function  
 \begin{equation}
\partial\Phi (\lambda,\mu)/\partial\lambda=   \frac{
  f(\mu) - f(\lambda) -f^\prime(\lambda)(  \mu -\lambda)
 } { (\lambda-\mu)^2 } .
  \label{eq:DO3a}\end{equation}
   is also bounded for $|\lambda|\leq R$, 
 $|\mu|\leq R$.
According to Lemma~\ref{Rmx}   the function
 \e{eq:DO4}  as well as   its derivatives in $\lambda$
are bounded in this region  provided $z=ia$ and $a$ is sufficiently large.

 It remains to consider the region
 $|\mu|\leq r$, $|\lambda|\geq R>r$ (and, similarly, $|\lambda|\leq r$, $|\mu|\geq R>r$).
 Here we use that   the denominator in  \eqref{eq:DO} is bounded from below because
    \begin{equation}
|g(\lambda) -g(\mu)|\geq (\mu^2+a^2)^{-m/2} - (\lambda^2+a^2)^{-m/2}
    \geq (r^2+a^2)^{-m/2} - (R^2+a^2)^{-m/2}\geq c>0.
\label{eq:DO5}\end{equation}
Therefore the function 
  \eqref{eq:DO} is   bounded. 
 Further, differentiating  \eqref{eq:DO}, we find that  
     \begin{equation}
\frac{\partial K(\lambda,\mu)}{ \partial\lambda}=   \frac{f^\prime(\lambda)}{g( \lambda) -g(\mu )}
  -  \frac{ ( f  (\lambda ) - f  (\mu ) ) g^\prime(\lambda)}{ (g( \lambda) -g(\mu ))^2 }.
\label{eq:DX}\end{equation}
Thus,  it follows from \eqref{eq:DO5}   for $|\mu|\leq r$,   $|\lambda|\geq R>r$ or 
for $|\lambda|\leq r$,   $|\mu|\geq R>r$ that
 \begin{equation}
\Bigl|\frac{\partial K(\lambda,\mu)}{ \partial\lambda} \Bigr| \leq  C(|f^\prime(\lambda)|+
( |f(\lambda)|+ |f(\mu)|)  |g^\prime(\lambda)|)  \leq  C_{1} (1+|\lambda|)^{-m-1}.
\label{eq:DX1}\end{equation}
This concludes the proof of Proposition~\ref{Rm1}.

\medskip

Now we check the estimates
\eqref{eq:DOI2} and \eqref{eq:DOI2x} under the assumptions of
    Proposition~\ref{Rm2}.  
    
    Let us first consider the region where $|\lambda| \leq 3r$, $|\mu|\leq 3r$
(with the square $|\lambda| \leq r$, $|\mu|\leq r$ removed).
Similarly to the proof of   Proposition~\ref{Rm1}, we have that functions
 \eqref{eq:DO3} and  \eqref{eq:DO3a} are bounded. The boundedness
 in this region of the function  \eqref{eq:DO4}, as well as of its derivatives in $\lambda$,
 follows, for sufficiently small $a$, from  Lemma~\ref{Rmx}.
 
 Next we consider the region $|\lambda| \geq 2r$, $|\mu|\geq 2r$
where  
\[
f(\lambda)-f(\mu)=(\lambda^m-i)^{-1}-(\mu^m-i)^{-1}=(\lambda^m-i)^{-1} (\mu^m-i)^{-1}(\mu-\lambda ) p(\lambda,\mu;0) 
\]
and the function $p(\lambda,\mu;0)$ is defined by formula \eqref{eq:DOI5}.
Hence   it follows from \eqref{eq:DO} and \eqref{eq:DO1} that   
$$
 K(\lambda,\mu)=\frac{(\lambda-z)^m}{\lambda^m-i}
\frac{(\mu-z)^m}{\mu^m-i}\frac{p(\lambda,\mu;0)}{p(\lambda,\mu;z)}.
$$
By virtue of estimate
\eqref{eq:DOI6a}, this function is bounded.
For the proof of \eqref{eq:DOI2x}, we use two estimates
\[
\Bigl| \frac{\partial}{\partial\lambda}\frac{(\lambda-z)^m}{\lambda^m-i}\Bigr|\leq C\lambda^{-2}
\]
and
\begin{equation}
\Bigl| \frac{\partial}{\partial\lambda}  \frac{p(\lambda,\mu;0)}{p(\lambda,\mu;z)}
\Bigr|\leq C(|\lambda|+|\mu|)^{-2}.
 \label{eq: DL4}\end{equation}
 The first of them is obvious. For the proof of the second, we remark that
\begin{eqnarray}
 \frac{\partial}{\partial\lambda}  \frac{p(\lambda,\mu;0)}{p(\lambda,\mu;z)}=
 p(\lambda,\mu;z)^{-2}
 \nonumber\\
  \Bigl( (p(\lambda,\mu;z)-p(\lambda,\mu;0))p^\prime_{\lambda}(\lambda,\mu;z)-(p^\prime_{\lambda}(\lambda,\mu;z)-p^\prime_{\lambda}(\lambda,\mu;0))p(\lambda,\mu;z)\Bigr).
 \label{eq: DL5} \end{eqnarray}
 Let us now take into account that $p(\lambda,\mu;z)$ is a polynomial in $\lambda$ and $\mu$ of degree $m-1$, but
  the terms of   order $m-1$ in 
 $p(\lambda,\mu;z)$ and $p(\lambda,\mu;0)$ cancel each other. Therefore their difference consists of terms $\lambda^p \mu^q$ where $p+q\leq m-2$. This gives the estimate
\begin{eqnarray}
 | p(\lambda,\mu;z)-p(\lambda,\mu;0)|\leq C (|\lambda|+|\mu|)^{m-2}.
 \label{eq: DL6} \end{eqnarray}  
Similarly, differentiating the difference $ p(\lambda,\mu;z)-p(\lambda,\mu;0)$, we see that 
\begin{eqnarray}
 \Bigl |  \partial  (p(\lambda,\mu;z)-p(\lambda,\mu;0))/\partial\lambda 
\Bigr |\leq C (|\lambda|+|\mu|)^{m-3}.
\label{eq: DL7} \end{eqnarray}
Substituting  \eqref{eq: DL6} and \eqref{eq: DL7} into  \eqref{eq: DL5}, we obtain
 \eqref{eq: DL4}.
 
 Let us finally consider the region where $|\lambda| \geq 3r$, $|\mu|\leq 2r$ or
$|\mu| \geq 3r$, $|\lambda |\leq 2r$. Here we use again representation 
\eqref{eq:DO} and  estimate 
\eqref{eq:DO5} which yields
\[
|K(\lambda,\mu)| \leq C (|f(\lambda)|+ |f(\mu)|)\leq C_{1}<\infty.
\]
The same bound \eqref{eq:DO5} shows that the derivative \eqref{eq:DX} satisfies estimate
\eqref{eq:DX1}, which yields \eqref{eq:DOI2x}.
 
 This concludes the proof of \eqref{eq:DOI2x} and hence that
 of Proposition~\ref{Rm2}.

\section{The Dirac operator}
 
 As an example to which Theorem~\ref{KRY}
directly applies, we
  now consider the Dirac operator describing a relativistic particle of
spin
$1/2$. Let  ${\cal H}=L_2({\Bbb R}^3; {\Bbb C}^4)$ and
\begin{equation}
 H_{00}=\sum_{j=1}^3 {\pmb \alpha}_j D_j +m {\pmb \alpha}_0,
\label{eq:3.3.11}\end{equation} 
where $m>0$ is the mass of  a particle and $4\times 4$ - Dirac
matrices satisfy the anticommutation relations
$$
  {\pmb \alpha}_i {\pmb \alpha}_j+{\pmb \alpha}_j {\pmb \alpha}_i=0,\quad i\neq j,\quad 
{\pmb \alpha}_i^2=I.
$$
 These relations determine the matrices ${\pmb \alpha}_i$ up to a unitary equivalence in the space
${\Bbb C}^4$. Their concrete choice is of no importance. If ${\pmb \alpha}_i$ are replaced by
$u {\pmb \alpha}_i u^\ast$ where
$u$ is  a unitary transformation in ${\Bbb C}^4$, then the operator $H_{00}$ is replaced by a unitary
equivalent operator of the same structure.

Making the Fourier transform $\Phi$, we find that $H_{00}=\Phi^\ast A \Phi$ where $A$ is multiplication by
the matrix  function (the symbol of $H_{00}$)
 $$A(\xi)=\sum_{j=1}^3 {\pmb \alpha}_j\xi_j +m {\pmb \alpha}_0.
   $$
 It is easy to see that
 $A(\xi)$ has
the eigenvalues
\begin{equation}
 a_{1,2}(\xi)=-a_{3,4}(\xi)=(|\xi|^2+m^2)^{1/2} 
\label{eq:3.3.eig}\end{equation}
of multiplicity $2$
so that 
\begin{equation}
A(\xi)= T(\xi)\Lambda (\xi) T^\ast (\xi)
\label{eq:3.3.14D}\end{equation}
 where the matrices $T(\xi)$ are unitary and
$\Lambda (\xi)=\mathrm{diag}\,\{a_1(\xi),a_2(\xi),a_3(\xi), a_4(\xi)\}$.
 A concrete form of the matrices $T(\xi)$ is inessential for us.
In particular,   the operator $H_{00}$  is selfadjoint on the Sobolev space
$ {\HH}^1({\Bbb R}^3; {\Bbb C}^4)=:{\cal D}(H_{00})$.

 Here we consider a couple of Dirac operators
  \begin{equation}
H_{0}=H_{00}+V_{0}, \quad H =H_{0}+V 
\label{eq:Di}\end{equation}
  where $V_{0}$ is multiplication by a symmetric bounded $4\times 4$ - matrix function $V_{0}(x)$ and a perturbation $V$ which is also a
  symmetric   $4\times 4$ - matrix function  satisfies the condition 
  \begin{equation}
 |V(x)|\leq C(1+|x|)^{-\rho},\quad \rho >3.
\label{eq:pot}\end{equation}
We do not make any special assumptions on matrices  $V_{0}(x)$ 
and $V (x)$. In particular, the spectrum of the operator $H_{0}$ might cover the whole real axis.
Particular cases 
\[
 V_{0}(x)=\sum_{j=0}^3v_j^{(0)}(x) {\pmb \alpha}_j + v^{(0)} (x),\quad  V(x)=\sum_{j=0}^3v_j(x)  {\pmb \alpha}_j + v (x),
\] 
where $ v_j^{(0)}, v ^{(0)}$ and  $v_j , v  $ are scalar functions, correspond   to an interaction of an electron (or of a
positron) with   magnetic and electric fields (with potentials  $(v_1,v_2, v_3)$ and $v $,
respectively). We suppose that   the ``background" potential $V_{0}(x)$ 
is only a bounded function  whereas     the perturbation $V (x)$ satisfies the condition    (\ref{eq:pot}). We shall show that in this case
 condition \eqref{eq:Rm} is satisfied for $m=3$.
 
 Differentiating the resolvent identity
$$
 R (z)-R_0 (z)=-R (z) V R_0 (z) ,\quad \Im z\neq 0 ,
$$
  we find that
\begin{equation}
 R^m(z)-R_0^m(z)= -\sum_{k=1}^m R^k(z)V  R_0^{m+1-k}(z).
\label{eq:4.2.5}\end{equation}    
If we replace  the resolvents $R $ and $R_0 $ by $R_{00} $ in the right-hand side of
\e{eq:4.2.5},
then we obtain the sum of terms $R_{00}^k  V  R_{00}^{m+1-k} $. It
is easy to see (this follows from Proposition~\ref{DIRAC} below) that, for $m\geq 3$,  these operators
  belong to the trace class. Therefore it suffices to justify the
replacement of  $R $ and $R_0 $ by $R_{00} $. However the boundedness of the operator $ 
(H_{00}-z)^k R ^k(z)$ is equivalent to the inclusion 
$$
  {\cal D}(H ^k)\subset {\cal D}(H_{00}^k). 
  $$
If $k>1$, this inclusion requires boundedness of derivatives of the function  $V $   and is, in
general, violated under the conditions above.  To bypass this difficulty, we suggest a trick based on commutation of the operators $\langle x \rangle^{-r} $ and $R(z)$. We introduce also the whole scale of symmetrically normed ideals ${\goth S}_{p}$.
  
  We start however  with operators $\langle x \rangle^{-r}  R^k _{00}(z)$.
  
  \begin{proposition}\label{DIRAC}
  If
  \begin{equation}
  p>3/\min\{r,k\}=:p(r,k), \quad p\geq 1,
  \label{eq:prk}\end{equation}
  then
    $$\langle x \rangle^{-r}  R^k _{00}(z)\in\mathfrak{S}_{p}.$$
     \end{proposition}
  
  \begin{proof} 
  Since $R^k _{00}(z)=\Phi^\ast (A-z)^{-k}\Phi$,
  we have to check that the integral operator with kernel
  $$ \langle x \rangle^{-r} \exp(i<x,\xi>)(A(\xi)-z)^{-k}$$
  belongs to the class ${\goth S}_{p}$ with $p$ determined by \e{eq:prk}.
  According to \e{eq:3.3.eig}, \e{eq:3.3.14D},
  the function 
  $$\langle \xi \rangle^k  (A(\xi)-z)^{-k}$$
  is bounded so that it suffices to consider the operator with kernel
   $$ \langle x \rangle^{-r} \exp(i<x,\xi>)\langle \xi \rangle^{-k}.$$ 
   This operator belongs to the required class ${\goth S}_{p}$ according, for example, to the results of \cite{RS}.
  \end{proof}

  Next we extend this result to the operator $H$ (or $H_{0}$). 
  
   \begin{proposition}\label{DIRACX}
If $p $ satisfies \e{eq:prk}, then
  $$\langle x \rangle^{-r}  R^k  (z)\in\mathfrak{S}_{p}.$$
 \end{proposition}
     
  \begin{proof}   
  If $k=1$, then
  $$\langle x \rangle^{-r}  R  (z)=(\langle x \rangle^{-r}  R_{00}  (z))\cdot
 ( (H_{00}-z)R  (z))$$
  and the second factor in the right-hand side is a bounded operator.
  Let us justify the passage from $k$ to $k+1$.
  Since 
  $$
  (H-z)  \langle x \rangle^{-r} - \langle x \rangle^{-r} (H-z)=
  [H_{00},\langle x \rangle^{-r}],
  $$
  we have that
  \begin{equation}
  \langle x \rangle^{-r}  R^{k+1}=
R  \langle x \rangle^{-r}  R^k + R [H_{00},\langle x \rangle^{-r} ] R^{k+1}. 
\label{eq:Dir}\end{equation}
 Let us write the first term in the right-hand side of \e{eq:Dir} as
  \begin{equation}
  R  \langle x \rangle^{-r}  R^k=(R  \langle x \rangle^{-r_{0}})( \langle x \rangle^{-kr_{0}}  R^k)
  \label{eq:Diri}\end{equation}
where $r_{0}=r(k+1)^{-1}$. Here $R  \langle x \rangle^{-r_{0}}\in
\mathfrak{S}_{p}$ for  $p>p(r_{0},1)$, and
$ \langle x \rangle^{-kr_{0}}  R^k \in
\mathfrak{S}_{p}$ for  $p>p(kr_{0},k)$. Therefore the product
\e{eq:Diri}  belongs to the class
$\mathfrak{S}_{p}$ where 
$$p^{-1}<p(r_{0},1)^{-1}+ p(kr_{0},k)^{-1}=p(r,k+1)^{-1}.$$
  The second term in the right-hand side of
\eqref{eq:Dir} is even better since the matrix $[H_{00},\langle x \rangle^{-r} ]$ is bounded by
 $\langle x \rangle^{-r-1}$.
  \end{proof}    

Now it easy to verify inclusion  \eqref{eq:Rm} for the Dirac operators.
      
\begin{theorem}\label{4.3.8d}
Suppose that the symmetric matrix functions $V_0 (x)$ and $V (x)$ are bounded and
$ V (x) $ satisfies condition  \eqref{eq:pot}  with $\rho >3$. 
Then condition \eqref{eq:Rm} is satisfied for all $m\geq 3$.
     \end{theorem}   

\begin{proof}
Indeed, let us write all terms in the right-hand side of \e{eq:4.2.5} as
\begin{equation}
  R^k V  R_0^{m+1-k} =(R^k  \langle x \rangle^{-r_{1}} )
( \langle x \rangle^\rho V) ( \langle x \rangle^{-r_{2}}R_0^{m+1-k})
\label{eq:4.vv}\end{equation} 
where $r_{1}=k\rho(m+1)^{-1}$, $r_{2}=(m+1-k)\rho(m+1)^{-1}$.
According to Proposition~\ref{DIRACX} the first term in the right-hand side of
\eqref{eq:4.vv} belongs to the class ${\goth S}_{p}$ for
$p>p_{1}=p(r_1,k)$ and the last term belongs to the class ${\goth S}_{p}$ for
$p>p_{2}=p(r_2, m+1-k)$. Since $p_{1} ^{-1}+   p_{2} ^{-1} >1$, the product 
\e{eq:4.vv} is trace class.
 \end{proof}

\begin{corollary}\label{4}
The  WO  $W_\pm (H,H_0)$    exist and are
complete.
 \end{corollary}

Finally, combining   Theorem~\ref{4.3.8d}  with Theorem~\ref{KRY}, we obtain 
 \begin{theorem}\label{KRYX}
 Let the operators $H_{0}$ and $H$ be given by formula \e{eq:Di} where 
 $H_{00}$ is the ``free" Dirac operator \e{eq:3.3.11}.
 Suppose that the symmetric matrix functions $V_0 (x)$ and $V (x)$ are bounded and
$ V (x) $ satisfies condition  \eqref{eq:pot}  with $\rho >3$. 
Let a  function $f(\lambda)$   have two bounded derivatives and
$$
   \partial^\alpha( f(\lambda)- f_{0}\lambda^{-3})=O(|\lambda |^{-3-\epsilon-\alpha}),
\quad \alpha=0,1,2, \quad \epsilon >0,
$$
where the constant $f_{0}$ is the same for $\lambda\rightarrow \infty$ and
$\lambda\rightarrow -\infty$.
 Then inclusion \eqref{eq:Rmf} holds and there exists a function $($the SSF$)$
 $\xi(\lambda;H,H_{0})$ satisfying condition \eqref{eq:xi2} where $m=3$  such that the trace formula \eqref{eq:TF}  is true.
   Moreover, for the  corresponding scattering matrix $S(\lambda;H,H_{0})$, the operator $S(\lambda;H,H_{0})-I\in\mathfrak{S}_{1}$ and relation
 \e{eq:Sxi} holds
for almost all $\lambda$ from the core of the spectrum of the operator $H_0$.
  \end{theorem}

      \end{document}